\documentclass[12pt]{amsart}

\def\@typesizes{%
       \or{5}{6.5}\or{6}{7.5}\or{7}{8.5}\or{8}{11}\or{9}{12}%
       \or{10}{13}
       \or{\@xipt}{14}\or{\@xiipt}{15}\or{\@xivpt}{18}%
       \or{\@xviipt}{20}\or{\@xxpt}{24}}

\usepackage{amsthm}
\usepackage{amsmath}
\usepackage{amssymb}
\usepackage{graphicx}
\usepackage{enumerate}
\usepackage{color}
\allowdisplaybreaks
\numberwithin{equation}{section}
\numberwithin{figure}{section}

\theoremstyle{plain}
\newtheorem{theorem}{ Theorem}[section]
\newtheorem{proposition}[theorem]{ Proposition}
\newtheorem{lemma}[theorem]{ Lemma}
\newtheorem{corollary}[theorem]{ Corollary}
\newtheorem{example}[theorem]{ Example}
\newtheorem{remark}[theorem]{ Remark}
\newtheorem{definition}[theorem]{ Definition}
\newtheorem{conjecture}{ Conjecture}

\def\BET{\begin{theorem}}
\def\ENT{\end{theorem}}
\def\BEP{\begin{proposition}}
\def\ENP{\end{proposition}}
\def\BEL{\begin{lemma}}
\def\ENL{\end{lemma}}
\def\BEC{\begin{corollary}}
\def\ENC{\end{corollary}}
\def\BEE{\begin{example} \rm}
\def\ENE{\end{example}}
\def\BER{\begin{remark} \rm}
\def\ENR{\end{remark}}
\def\BED{\begin{definition} \rm}
\def\END{\end{definition}}
\def\BECJ{\begin{conjecture}}
\def\ENCJ{\end{conjecture}}

\def\bea{\begin{eqnarray}}
\def\eea{\end{eqnarray}}

\def\beas{\begin{eqnarray*}}
\def\eeas{\end{eqnarray*}}

\def\beq{\begin{equation}}
\def\eeq{\end{equation}}

\def\beal{\begin{align*}}

\def\eeal{ \end{align*} }

\def\roweq{\nonumber \\ &=& }
\def\rowleq{\nonumber \\  & \leq & }

\def\rowpl{\nonumber \\ & + & }
\def\rowmi{\nonumber \\ & - & }

\def\bfV{{\bf V}}

\def\bbC{{\mathbb C}}

\def\bbN{{\mathbb N}}

\def\bbR{{\mathbb R}}

\def\cA{{\mathcal A}}

\def\cD{{\mathcal D}}

\def\cH{{\mathcal H}}

\def\cK{{\mathcal K}}

\def\cO{{\mathcal O}}

\def\cT{{\mathcal T}}

\def\cZ{{\mathcal Z}}

\def\ef{\eqref}
\def\sfk{{\sf k}}

\begin{document}

\title[Plummeting eigenvalues]{Plummeting and blinking eigenvalues of the Robin Laplacian in a cuspidal domain}



\author{Sergei A. Nazarov}
\address{
Saint-Petersburg State University,
Universitetskaya nab., 7--9,  St. Petersburg, 199034, Russia, and  \hfill\break
Institute of Problems of Mechanical Engineering RAS,
V.O., Bolshoj pr., 61, St. Petersburg, 199178, Russia}
\email{s.nazarov@spbu.ru, srgnazarov@yahoo.co.uk}

\author{Nicolas Popoff}
\address{Institut de Math\'ematiques de Bordeaux UMR 5251,
Universit\'e de Bordeaux, 351 cours de la Lib\'eration - F 33 405 Talence, France} 
\email{nicolas.popoff@u-bordeaux.fr}

\author{Jari Taskinen}
\address{Department of Mathematics and Statistics, P.O.Box 68, 
University of Helsinki, 00014 Helsinki, Finland}
\email{jari.taskinen@helsinki.fi}

\thanks{ The first  named author was supported by the grant 
17-11-01003 of the Russian Science Foundation.} 

\begin{abstract}We consider the Robin Laplacian in the domains
$\Omega$ and $\Omega^\varepsilon$, $\varepsilon >0$,  with sharp and blunted 
cusps, respectively.  Assuming that the Robin coefficient $a$ is large enough,
the spectrum of the problem in $\Omega$ is known to be residual and to
cover the whole complex plane, but on the contrary, the spectrum in the 
Lipschitz domain $\Omega^\varepsilon$ is discrete. However, our results reveal 
the strange behavior of the discrete spectrum as the blunting parameter 
$\varepsilon$ 
tends to 0: we construct asymptotic forms of the eigenvalues and detect  families of "hardly movable" and "plummeting" ones. The first type
of the eigenvalues do not leave
a small neighborhood of a point for any small $\varepsilon > 0$ while the 
second ones move at a high rate $O(|\ln \varepsilon|)$ downwards along the
real axis $\bbR$ to $ -\infty$. At the same time, any point $\lambda \in \bbR$
is a "blinking eigenvalue", i.e.,  it belongs to the spectrum of the problem
in $\Omega^\varepsilon$ almost periodically in the $|\ln \varepsilon|$-scale.
Besides standard spectral theory, we use the techniques of dimension 
reduction and  self-adjoint extensions to obtain these results.
  
\end{abstract}
\maketitle

\section{Introduction.}
\label{sec1}
\subsection{Formulation of the problems.}
\label{sec1.1}

\begin{figure}
\begin{center}

\includegraphics[width=10cm]{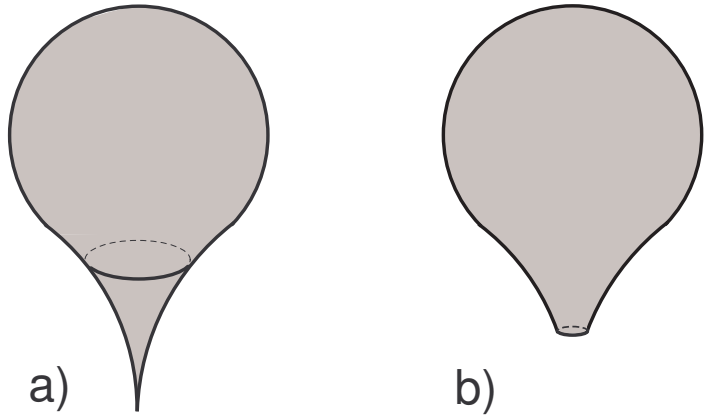}
\end{center}
\caption{Problem domain with a cusp (a) and domain with a blunted csup (b)}\label{fig1}
\end{figure}

We consider a family of spectral problems for the Laplace operator with the Robin 
condition
\bea
- \Delta u^\varepsilon (x)&=& \lambda^\varepsilon u^\varepsilon(x), \ \   x \in \Omega^\varepsilon,
\label{1} \\ 
\partial_\nu u^\varepsilon (x) &=& a u^\varepsilon(x), \ \ \  x \in \partial\Omega^\varepsilon,
\label{2}
\eea
in the domain (Fig.\,\ref{fig1},b)
\bea
\Omega^\varepsilon = \Omega \smallsetminus \overline{\Pi^\varepsilon}
\subset \bbR^n \ , \ \ n\geq 2 , \label{3}
\eea
where $\varepsilon \in (0,\varepsilon_0]$ is a small parameter, $a\in \bbR$ a constant, $\lambda^{\varepsilon}$ the spectral parameter, 
$\partial_\nu$ is the outward normal derivative,
\bea
\Pi^d = \{ x = (y,z) \in \bbR^{n-1} \times \bbR \, : \, z = x_n \in (0,d) 
, \eta = z^{-2 }y \in \omega\} \ ,  \ \ d > 0 , \label{4}
\eea
$\omega$ is a domain in $\bbR^{n-1}$ with Lipschitz boundary and compact closure
$\overline{\omega} = \omega \cup \partial \omega$, and $\Omega$ is assumed to 
coincide with the cusp $\Pi^d$ in a neighborhood of the coordinate origin $\cO$ (Fig.\,\ref{fig1},a). The domain $\Omega$ is Lipschitz everywhere, except at 
the point $\cO$.

For $\varepsilon > 0$ the domain \eqref{3} is Lipschitz and the spectrum of the 
problem \eqref{1}--\eqref{2} is discrete, consisting of the monotone increasing 
unbounded sequence of eigenvalues
\bea
\lambda_1^\varepsilon < \lambda_2^\varepsilon 
\leq \lambda_3^\varepsilon \leq \ldots \leq \lambda_m^\varepsilon \leq \ldots
\to + \infty . \label{5}
\eea
As studied for example in \cite{Dan00}, it is possible to define a limit problem ($\varepsilon =0$) 
in the cuspidal domain $\Omega = \Omega^0$, 
\bea
- \Delta u(x)&=& \lambda u(x), \ x \in \Omega \ \ , \ \ \ 
\partial_\nu u (x) = a u(x), \ x \in \partial\Omega ,
\label{6}
\eea
Moreover, it is known that if the constant coefficient in \eqref{2} 
is non-positive, this problem has dicrete spectrum. When $a$ is positive, it was proven in \cite{na549} that the discrete spectrum
constitutes the whole spectrum $\sigma$ of \eqref{6}  only in the case  $a < a_\dagger$, 
while $\sigma$ becomes the residual spectrum and covers the whole complex plane $\bbC$ in the case
\bea
a \geq a_\dagger = \Big( n -\frac32 \Big)^2 \frac{|\omega|}{|\partial \omega|} ,
\label{7}
\eea
where $|\omega| = \mbox{mes}_{n-1} \omega$ is the volume of the cross-section
and  $|\partial\omega| = \mbox{mes}_{n-2} \partial \omega$ is the area of its boundary.

\subsection{State of art.}
\label{sec1.2}

When the parameter $a$ is positive, it is not clear how to present a reasonnable weak formulation of the limit problem ($\varepsilon=0$) in the cusp. The Robin Laplacian has been studied in arbitrary domains in different ways. Using a variational approach, a possible start is to consider the quadratic form
$$ 
u \mapsto \|\nabla u; L^{2}(\Omega) \|^2-\|a^{1/2} u, L^{2}(\partial\Omega)\|^2,
$$
defined on its natural domain, as shown in \cite[Section 3]{ArWa03} and \cite{Dan13}. In these works, measures on the boundary are considered, including our case $a>0$. For a domain with a cusp, the resulting operator is not necessarily self-adjoint. If the cusp is, roughly, {\it less sharp than quadratic}, then the form is bounded from below, and the spectrum is discrete (see \cite{MaPo,na478,na549} and, e.g., \cite{KovPan18} for a recent study 
of the corresponding eigenvalue sequence itself). But, as it was shown
in \cite{na478,na549}, the nature of the problem operator may become 
completely different, as it may lose its semi-boundedness, if the cusp is {\it sharper than quadratic}, see also \cite[Section 5]{Dan13}. For the critical case of a quadratic cusp \eqref{4} considered here, the spectrum is discrete if and only if $a<a_{\dagger}$, since the spectrum becomes residual and fills in the 
whole complex plane when $a \geq a_{\dagger}$, see \cite{na549}.

Let us review the {\it Steklov} problem related to \eqref{1}--\eqref{2}.
In the paper \cite{na401} it was shown that the spectrum (subset of $\bbR_+$)  of the 
Steklov problem  in a domain with a peak type boundary singularity is either discrete or
may contain a continuous component depending on the sharpness of the peak. 
Related to this, the linear water wave problem, which contains the Steklov 
condition on a part of the boundary, was considered  in \cite{na493} in 
domains with  rotational cusps: a formulation of the problem as a
Fredholm operator of index zero was given with the help of appropriate
radiation conditions, and it was proven  that the continuous 
spectrum is non-empty and consist of the ray $[\lambda_\dagger ,+ \infty) 
\subset \bbR$ with a certain cut-off point $\lambda_\dagger \geq 0$.

The reference \cite{NaTaNEW} contains a study of the Laplace equation 
\bea
- \Delta u^\varepsilon (x) = 0, \ \ x \in \Omega^\varepsilon , \label{SL1}
\eea
with the spectral Steklov and Dirichlet boundary condition
\bea
\partial_\nu u^\varepsilon(x) &=& \lambda^\varepsilon u^\varepsilon (x) \ ,
\ \ x \in \partial \Omega^\varepsilon \smallsetminus \overline{\omega^\varepsilon} ,
\nonumber \\
u^\varepsilon(x) &=& 0  \ , 
\ \ x \in \omega^\varepsilon , \label{SL2}
\eea
where $\omega^\varepsilon = \{ x \in \Pi^d \, : \, x = \varepsilon \}$ is the
end of the blunted cusp. The spectrum of this problem is discrete and similar to
\eqref{5}. According to \cite{na401}, the spectrum of the limit Steklov problem 
($\varepsilon =0$) in the cuspidal domain $\Omega$ is continuous and equals 
$[\lambda_\dagger , + \infty)$ with the cut-off value $\lambda_\dagger = a_\dagger$,
\eqref{7}, while it was discovered in \cite{NaTaNEW} that the eigenvalues 
$\lambda_m^\varepsilon > \lambda_\dagger$ of the Steklov-Dirichlet in
$\Omega^\varepsilon$ behave "strangely" as $\varepsilon \to +0$, namely they
"glide" within the semi-axis $(\lambda_\dagger , + \infty)$ at a 
high rate $O(|\ln \varepsilon|)$, which however slows down near 
$\lambda_\dagger$ so as to make $\lambda_m^\varepsilon$ "parachute" smoothly 
on $\lambda_\dagger$. Moreover, each point $\lambda > \lambda_\dagger$
constitutes a "blinking" eigenvalue of the problem \eqref{SL1}, \eqref{SL2}, namely,
for every $\lambda > \lambda_\dagger$ there exists a positive sequence
$\{ \varepsilon_k(\lambda) \}_{k=1}^\infty$ tending to 0 such that 
$\lambda$ becomes a true eigenvalue for  the  problem \eqref{SL1}, 
\eqref{SL2} in the 
domain $\Omega^{\varepsilon}$ for some $\varepsilon$  close to
$\varepsilon_k(\lambda)$, for any $k$. This phenomenon can be used  to 
construct a singular Weyl sequence at $\lambda$ 
for the Steklov problem operator in $\Omega$, which provides  a novel mechanism  to
form the continuous spectrum from a family of discrete spectra.

\subsection{Outline of the paper.}
\label{sec1.3}
The asymptotic expansions of the solutions of the problem \eqref{6}  near the tip 
$\cO$  were derived in detail in \cite{na549} and will be reproduced in Sections 
\ref{sec2.1}--\ref{sec2.3}, and although they are the same as in the case of the
Steklov problem in  \cite{NaTaNEW}, the rest of the material is quite different.
In Section \ref{sec2.4} we determine  all self-adjoint extensions of the Robin-
Laplacian, which is originally defined in  the small domain \eqref{31}. However, 
none of the extensions is lower semi-bounded (for general results 
on non-semi-bounded sesquilinear, see \cite{DFHW, MI}),
which somehow reflects the fact that the spectrum of the problem \eqref{6} covers the whole
complex plane $\bbC$, see \cite{na549}.  Some of these extensions 
$A_{\theta_\bullet^\varepsilon}$ have a peculiar property, namely their eigenfunctions
leave a relatively small discrepancy in the Robin condition at the end of the blunted cusp,
see Section \ref{sec4.1}, and thus can be regarded as good candidates
to model the singularly perturbed problem \eqref{1}, \eqref{2}, cf. the argumentation
in \cite{na240}. This plan will be realized in Sections \ref{sec4.2}--\ref{sec4.4}, 
where it is shown that a small neighborhood of any point of the spectrum of the 
extension $A_{\theta_\bullet^\varepsilon}$  contains an eigenvalue of the problem 
\eqref{1}, \eqref{2}  in $\Omega^\varepsilon$.

An important property is that the extension parameter $\theta_\bullet^\varepsilon$ 
in \eqref{63} is a periodic function in the logarithmic scale $|\ln \varepsilon|$, 
hence,  the spectrum $\sigma ( A_{\theta_\bullet^\varepsilon} )$ gains the 
same property. Among the eigenvalues of $\sigma( A_{\theta_\bullet^\varepsilon})$ 
there exist the so called stable eigenvalues, which are hardly movable and are generated
by "trapped modes", i.e., solutions of the problem \eqref{6} in the Sobolev
space $H^1(\Omega)$. However, according to Theorem \ref{T3}
there certainly exist also eigenvalues of  $\sigma( A_{\theta_\bullet^\varepsilon})$
which are generated by "diffraction" solutions \eqref{46} of the problem \eqref{6},
move downwards at a high speed along the real axis according to the formulas \eqref{58} 
and \eqref{64}, and therefore are called "plummeting". In other words,
the spectrum  $\sigma( A_{\theta_\bullet^\varepsilon})$ is indeed periodic in
$|\ln \varepsilon|$, although as a set only, because some points of it 
move purposefully to a fixed direction as $\varepsilon \to +0$. Such a situation
may occur only in a situation, when the model operator is not lower semi-bounded.
This does not happen in the case of the Steklov problem, which is investigated in
\cite{NaTaNEW} and characterized by the phenomenon of "gliding" eigenvalues
(see the end of Section \ref{sec1.2}).

\section{Theorem on asymptotics in the cuspidal domain.}
\label{sec2}
\subsection{Formal asymptotics.}
\label{sec2.1}

We aim  to  present the asymptotics of solutions of the problem
\bea
- \Delta u - \lambda u(x) &=& f(x) , \ x \in \Omega,
\label{11} \\ 
\partial_\nu u-  a u(x) &=& 0 , \ x \in \partial\Omega,
\label{12}
\eea
where $a$ satisfies \eqref{7}, and first of all we will describe a formal 
procedure under the assumptions that the boundary $\partial \omega$ is 
smooth and the right-hand side $f$ vanishes near the cuspidal tip $\cO$.

Since the diameter $O(\zeta^2)$ of the cross-section
\beas
\omega^\zeta = \{ (y,z) \in \Pi^d \, : \, z = \zeta \} \ , \ \ 
0 < \zeta \ll d \  , 
\eeas
is much less than its distance $\zeta$ to the tip $\cO$,  it is logical to 
accept the standard asymptotic ansatz in thin domains,  see e.g. \cite[Ch.\,14]{MaNaPl},
\bea
u(x) = w(z) + W(\eta, z) + \ldots , \label{13}
\eea
where the dots stand for inessential higher order terms, $\eta = z^{-2} y$ is the 
"rapid" variable used in \eqref{4}, and the power-law functions
\bea
w(z) = z^\mu w_0 \ , \ \ W(\eta,z)  = z^{\mu+2} W_0(\eta) \label{14},
\eea
where $W_{0}\in H^{1}(\omega)$, are to be determined. We insert the  ansatz
\eqref{13}  into the differential equation \eqref{11}, extract terms of order 
$z^{\mu-2}$ as $z \to +0$, and thus obtain the relation
\bea
- \Delta_\eta W(\eta, z) = z^4 \partial_z^2 w(z) \ , \ \ \eta \in \omega. \label{15}
\eea
The normal derivative on the lateral side 
$\Gamma^d =  \{ x \, : \,  \eta \in \partial \omega, \ z \in (0,d) \}$ 
of the cusp $\Pi^d$ equals

\beas
\partial_\nu = \big( 1 + 4 z^2 |\eta \cdot \nu'(\eta)|^2 \big)^{-1/2} 
\big( z^{-2} \nu'(\eta) \cdot \nabla_\eta - 2 z\eta\cdot\nu ' (\eta)  \partial_z+2(\eta\cdot\nu'(\eta))(\eta\cdot\nabla_{\eta}) \big) ,
\eeas
where $\nu'(\eta)$ is the unit normal vector on the boundary of the domain $\omega \subset \bbR^{n-1}$ and the central dot stands for the scalar product in the Euclidean spaces.
Hence, by considering the order $z^\mu$ in  the boundary condition \eqref{12},  we derive
the relation 
\bea
\nu'(\eta) \cdot \nabla_{\eta} W(\eta, z) = 2 \eta \cdot \nu'(z) z^{3} \partial_z w(z) 
+a z^{2 }w(z) \ , \ \ \eta \in \partial \omega. \label{16}
\eea
Using the formula 
\beas
\int\limits_{\partial \omega} \eta \cdot \nu'(\eta) ds_\eta = \int\limits_\omega
\nabla_\eta  \cdot \eta \, d\eta =   (n-1) |\omega| 
\eeas
we see that the compatibility condition in the Neumann problem \eqref{15}, \eqref{16}
reads as
\beas
0 &=& \int\limits_\omega \partial_z^2 w(z) d \eta + \int\limits_{\partial\omega} \big( 2 \eta 
\cdot \nu'(\eta) z^3 \partial_z w(z) + a z^2 w(z) \big) ds_\eta
\roweq
|\omega| \partial_z^2  w(z) + 2(n-1) |\omega| z^{3} \partial_z w(z) + a z^ {2} 
|\partial\omega| w(z) 
\eeas
and turns into the ordinary differential equation of Euler type 
\bea
- \frac{d}{dz} \Big( z^{2(n-1)} \frac{dw}{dz} (z) \Big) 
= A z^{2(n-2)} w(z) \ , \ \ z >0,
\label{17}
\eea
where
\beas
A = a\frac{|\partial \omega|}{|\omega|}. 
\eeas

The general solution of the equation \eqref{17} is of the form
\bea
w(z) = b_+ w^+(z) + b_- w^- (z)\ , \ \ b_\pm \in \bbC , \label{19}
\eea
where we have in the case  $a > a_\dagger$, \eqref{7},
\bea
w^\pm (z) = w_0 z^{\pm i \mu_0 - n + 3/2} \ \ \mbox{with} \ 
\mu_0 = \sqrt{ A - \Big(n -\frac32 \Big)^2 } > 0 \ , \label{20}
\eea
and in the case  $a = a_\dagger$, $A=(n-3/2)^2$
\bea
w^\pm (z) = w_0 z^{-n +3/2} \big( \ln z \mp i\big). \label{21}
\eea
The normalization factor $w_0$ and the peculiar form of \eqref{21}
will be clarified later on. 

Since the compatibility condition in the problem \ef{15}--\ef{16} is fulfilled, there
exists a solution $W$ defined up to an additive constant with respect to $\eta$. 
To make the solution unique, we impose the orthogonality condition 
\bea
\int\limits_\omega W(\eta,z)d\eta = 0. \label{22}
\eea

\subsection{Weak formulation of the problem.}
\label{sec2.2}
We introduce the weighted Sobolev space $V_\beta^1(\Omega)$ as the completion of 
the linear space $C_c^\infty (\overline \Omega \smallsetminus \cO )$ (infinitely 
differentiable functions vanishing in a neighborhood of the point $\cO$) with 
respect to the norm
\bea
\Vert u ; V_\beta^1(\Omega) \Vert = \big( 
\Vert r^\beta \nabla u ; L^2(\Omega) \Vert^2  + \Vert r^{\beta-1} u ; L^2(\Omega) \Vert^2
\big)^{1/2} \label{23}
\eea
where $r = {\rm dist}\,(x, \cO)$ and $\beta \in \bbR$ is a weight index. The weighted
Lebesgue space $V_\beta^0 (\Omega)$ is endowed with the norm 
$\Vert r^\beta  u ; L^2(\Omega) \Vert$.

\BER
\label{R1}
The norm \ef{23} is the same as the classical Kondratiev norm \cite{Ko}, but the 
reason for the  use of this norm in \cite{na549} as well as in the present paper is 
not the conventional one, since  the shape of the domain $\Omega$ near the singularity point $\cO$ 
is not conical nor angular as in Kondratiev's works.  This can be seen for example in the asymptotic ansatz for 
solutions: $W_{0}$ being in $H^{1}(\omega)$, the sum 
$w(z) + W(z^{-2}y,z)$,  see \ef{14}, belongs to $V_\beta^1(\Pi^d)$, if and only if 
\beas
\beta > - {\rm Re}\,\mu - n + 3/2
\eeas
even in the case $w =0$. However, if $w=0$ and $W_0(\eta)$ is independent of
the fast variable $\eta = z^{-2} y $, the condition for the space $V_\beta^1 (\Pi^d)$ 
to include $W$ becomes much less restrictive:
\beas
\beta > - {\rm Re}\,\mu - n - 1/2 . \hfill \hskip2cm  \boxtimes
\eeas 
\ENR

According to \cite{na549}, the weak formulation of the problem \ef{11}-\ef{12}
for the unknown $u \in V_\beta^1(\Omega)$ reads as the integral identity
\bea
(\nabla u, \nabla v)_\Omega - \lambda (u,v)_\Omega - a(u,v)_\Omega
= f(v) \ \ \ \forall\, v \in V_{-\beta}^1(\Omega) , \label{24}
\eea
where $f \in V_{-\beta}^1(\Omega)^* $ is an (anti)linear functional on 
$V_{-\beta}^1 (\Omega)$, in particular
\bea
f(v) = (f,v)_\Omega \ \ \ \mbox{with} \ f \in V_{\beta+1}^0 (\Omega). \label{25}
\eea
Here $( \cdot , \cdot )_\Omega$ is the natural scalar product in $L^2(\Omega)$, extended
by density to the duality between the spaces $V_\beta^0(\Omega)$ and $V_{-\beta}^0(\Omega)$.
According to definition \ef{23} and the weighted trace inequality 
\cite[Lemma\,2.2]{na549}
\beas
\Vert r^\beta u ; L^2(\partial\Omega)\Vert \leq c \Vert u; V_\beta^1(\Omega)\Vert,
\eeas
all expressions in the integral identity \ef{24} are properly defined so that it 
determines a continuous mapping 
\beas
V_\beta^1(\Omega) \ni u \mapsto \cT_\beta(\lambda) u = 
f \in V_{-\beta}^1(\Omega)^* . 
\eeas
We observe that for every $\beta$,  $ \cT_{-\beta}(\lambda)$ is the adjoint operator of $ \cT_\beta(\lambda)$. 
In Section \ref{sec3} we use the arguments of \cite{na549} to  describe the properties of $ 
\cT_{\pm \beta}(\lambda)$ in the particular case $\beta =1$.

\subsection{Theorem on asymptotics.}
\label{sec2.3}
We consider the problem \ef{11}--\ef{12} with the right-hand side
\bea
f \in L^2(\Omega) \subset V_{-1}^1 (\Omega)^*  \ \ 
\label{28}
\eea
(i.e. $\beta= -1$ in \ef{25}) and its solution $u \in V_1^1(\Omega) \subset L^2(\Omega) . $ 

The following assertion was verified in \cite{na549}. 

\BET\label{T1}
If \ef{7} and \ef{28} hold true, the above mentioned solution has 
the asymptotic form
\bea
u(x) = \chi(x) \big( w(z) + W(z^{-2} y,z) \big) + \widetilde u(x) , 
\label{29}
\eea
where $\chi$ is a smooth cut-off function which is equal to $1$  in $\Pi^{d/2}$
and $0$ in $\Omega^d$, see \ef{4} and \ef{3}. The term $w$ in \ef{29} is the linear
combination \ef{19} with some coefficients $b_\pm$ and terms \ef{20} or \ef{21}, and 
$W\in H^{1}(\Omega)$ is the solution of the problem \ef{15}, \ef{16}, \ef{22}. The coefficients 
$b_\pm$ and the remainder $\widetilde u \in V_{-1}^1 (\Omega)$ 
satisfy the estimate
\bea
\big( |b_+|^2 + |b_-|^2 + \Vert \widetilde u ; V_{-1}^1 (\Omega) \Vert^2 \big)^{1/2}
\leq c\big( \Vert f ; L^2(\Omega) \Vert + \Vert u; V_1^1(\Omega)\Vert \big) ,
\label{30}
\eea
here the factor  $c>0 $ is independent of $f$ and $u$. 
\ENT

\BER
\label{R2}
According to formulas \ef{20}, \ef{21} and Remark \ref{R1}, the detached asymptotic term
on the right of \ef{29} belongs to the space $V_\gamma^1(\Omega)$ with any $\gamma > 0$, but 
it is not contained in $V_{-1}^1(\Omega)$. Furthermore, as for second derivatives we have
$\nabla^2 w \in V_{\gamma+1}^0(\Pi^d) $ and $\nabla^2 W \in V_{\gamma+2}^0(\Pi^d) $, but in
general $\nabla^2 W \notin V_{\gamma+1}^0(\Pi^d) $. As it was verified in \cite{na549},  
for the solution  $u$ 
there holds  $\nabla^2 u \in V_{1+2}^1(\Omega)$  and $\nabla^2 \widetilde u \in 
V_{-1+2}^1(\Omega)$. We emphasize that the term
$W$, generated  by \ef{20}, \ef{19} and \ef{16} is defined up to the addendum
\beas
b_+^0 z ^{+i \mu -n + 7/2} + b_-^0 z ^{-i \mu -n + 7/2} 
\eeas
which is independent of $\eta$ and belongs to $V_1^1(\Omega)$ and can therefore be omitted 
in the asymptotic representation \ef{29}; this was the very reason for imposing the orthogonality condition \ef{22}. 

All these peculiarities again underline the difference of the conical \cite{Ko} and
cuspidal \cite{na549} irregularities of boundaries.
\ENR

By $\bfV_{\pm 1}^1(\Omega)$ we denote the weighted  space with detached asymptotics
(see \cite[Ch.\,6]{NaPl}, \cite[Sect.\,3]{na262} and others),
which consists of functions of the form \ef{29} and endow it with the norm on the left of
\ef{30}; the Hilbertian structure of this norm  can also be identified with the 
direct product
\bea
\bbC^2 \times V_{- 1}^1(\Omega) \ni (b_\pm, \widetilde u) , \label{c2v}
\eea 
although these will not be used later on.

\subsection{Symmetric and self-adjoint operators.}
\label{sec2.4}
As in \cite{na493} we associate to the problem \ef{11}--\ef{12} the symmetric operator
$ \cA$ in $L^2(\Omega)$, which has the differential expression $- \Delta$ and
the domain
\bea
\cD( \cA) = \{ u \in V_{-1}^1(\Omega) \, : \, \Delta u \in L^2(\Omega), \ \partial_n u
= au \ \mbox{on} \ \partial \Omega \smallsetminus \cO \}.
\label{31}
\eea
Notice that the inclusions in \eqref{31} assure that $u \in H^2(\overline \Omega
\smallsetminus \cO)$ and therefore the trace of $\partial_n u$ is properly 
defined on $\partial \Omega \smallsetminus \cO$. 
By Theorem \ref{T1}, see also \cite[Prop.\,3.11]{na549}, the adjoint operator $ \cA^*$
has the same differential expression but a larger domain
\beas
\cD( \cA^* ) = \{ u \in V_{1}^1(\Omega) \, : \, \Delta u \in L^2(\Omega), \ \partial_n u
= au \ \mbox{on} \ \partial \Omega \smallsetminus \cO \}.
\eeas
In view of Theorem \ref{T1} on asymptotics, the dimension of the  quotient space $\cD( \cA^*) / \cD( \cA)$ equals 2.

In order to describe the self-adjoint extensions of the operator $ \cA$ we reproduce
a calculation from \cite[\S3.4]{na549}. Let $\zeta > 0$. Applying the Green formula in the domain $\Omega^\zeta$ for  $u^1$, $u^2 \in \cD( \cA^*)$ and sending $\zeta$ to $+0$, we get
\bea
& & ( \cA^* u^1,u^2)_\Xi -  (u^1,  \cA^*  u^2)_\Xi
\roweq
\lim\limits_{\zeta \to +0} \int\limits_{\omega^\zeta} \big( \overline{u^2(y,z)}
\partial_z u^1(y,z) - u^1(y,z) \overline{\partial_z u^2(y,z)} \big)
\Big|_{z = \zeta} dy . \label{32}
\eea
Substituting $V_\beta^1(\Omega) \ni u^j \mapsto \widetilde u^j  \in 
V_{-\beta}^1(\Omega) $ at least for one of the indices $j=1,2$ makes the limit on the right
of \ef{32} equal to zero. Hence, we can replace $u^j$ by $w^j + W^j$ in the representation 
formula \ef{29}. Furthermore, since $W^j$ has the additional factor $z^2$, cf. \ef{14},
we can neglect the second term in this sum and write 
\beas
& & ( \cA^* u^1,u^2)_\Xi -  (u^1,  \cA^*  u^2)_\Xi
\roweq
\lim\limits_{\zeta \to +0} \zeta^{2(n-1)} 
\int\limits_{\omega} \big( \overline{w^2(\zeta)}
\partial_z w^1(\zeta) - w^1(\zeta) \overline{\partial_z w^2(\zeta)} \big)
 d\eta . 
\eeas
Finally, we fix the normalization factor in \ef{20} and \ef{21}, 
\begin{equation*}
w_0 = \left\{
\begin{array}{ll}
\sqrt{2 \mu_0 |\omega|} \ \ & \mbox{for} \ a > a_\dagger , \\
\sqrt{2  |\omega|} \ \ & \mbox{for} \ a = a_\dagger ,
\end{array}
\right. 
\end{equation*}
and  obtain 
\bea
& & ( \cA^* u^1,u^2)_\Xi -  (u^1,  \cA^*  u^2)_\Xi = i \big( \overline{b_+^2} b_+^1
- \overline{b_-^2} b_-^1 \big) , \label{34}
\eea
where $b_\pm^j$ are the coefficients of the linear combination $w_j$,
see \ef{19}. Repeating a traditional argument in \cite{RoBe} we observe that if $u^1$ and $u^2$ 
belong to the domain of a self-adjoint operator, then the left-hand side of \ef{34} vanishes 
and thus the coefficients must be related as $b_-^j = e^{i \theta} b_+^j$ for 
some $\theta \in [0, 2 \pi) $. 

\BET
\label{T2}
The restriction $ \cA_\theta$, where $\theta \in [0, 2 \pi)$, 
of the operator $ \cA^*$ to the subspace 
\bea
\big\{ u \in \cD( \cA^*) \, : \, b_+ = e^{i \theta} b_- \big\}  \label{35}
\eea
is a self-adjoint extension of the operator $ \cA$. Moreover, the domain of any self-adjoint 
extension of $ \cA$ equals \ef{35} for some parameter $\theta \in [0, 2 \pi)$. 
\ENT 

%
%
%

\section{Spectra of self-adjoint extensions.}
\label{sec3}
\subsection{Operator kernels.}
\label{sec3.1}
We fix the parameter $\lambda$, assume (1.7), and compare $\lambda$ and compare the Fredholm operators $ \cT_{+1} (\lambda)$ and 
$ \cT_{-1}(\lambda)$, which are adjoint to each other and therefore
\bea
{\rm dim \, ker} \, \cT_{\pm 1}(\lambda) = {\rm dim \, coker} \,  \cT_{\mp 1 }(\lambda)
\ \ \Rightarrow \ \  
{\rm Ind} \, \cT_{+  1}(\lambda) =  - {\rm Ind}\, \cT_{- 1}(\lambda) . \label{41}
\eea
Clearly, $V_{-1}^1(\Omega) \subset V_{+1}^1(\Omega)$ and
\bea
{\rm ker} \, \cT_{-1}(\lambda) \subset   {\rm ker}\, \cT_{+ 1}(\lambda) . \label{42}
\eea
Furthermore, Theorem \ref{T1} on asymptotics shows that 
\bea
{\rm Ind} \, \cT_{+  1}(\lambda) = 2+  {\rm Ind}\, \cT_{- 1}(\lambda)  \label{43}
\eea
where 2 is nothing but the number of the detached terms in formula \ef{29}, see \ef{19} with 
free constants $b_\pm$. From \ef{41} and \ef{43} we deduce that Ind\,$ \cT_{\pm} (\lambda) = \pm1$, and  taking \ef{42} into account yields
\beas
{\rm ker} \, \cT_{+1}(\lambda) =  {\rm ker}\, \cT_{- 1}(\lambda) \oplus \cZ , \ \ \ 
{\rm dim}\, \cZ = 1. 
\eeas
Any non-zero function $Z \in \cZ = {\rm ker}\, \cT_{+ 1} (\lambda) \ominus 
{\rm ker}\, \cT_{- 1} (\lambda)$, i.e. a solution of the homogeneous problem \ef{11}--\ef{12}
belonging to $V_1^1(\Omega)$, has the representation \ef{29} with the linear combination 
\ef{19}; the generalized Green formula \ef{34} with $u^1=u^2 = Z$
yields the equality 
\bea
0 = i |b_+|^2 - i |b_-|^2. 
\label{45}
\eea
If $b_\pm =0$, we arrive at the contradiction
\beas
Z \in V_{-1}^1(\Omega) \ \ \Rightarrow \ \ Z \in {\rm ker}\, \cT_{-1}(\lambda). 
\eeas
Thus, none of the coefficients vanishes and in view of \ef{45} we can choose a particular 
solution 
\bea
Z_\lambda(x) &=& \chi(x) \big( w_-(z) + W_-(z^{-2} y ,z)+e^{i \Theta (\lambda) } \big( w_+(z) + W_+(z^{-2} y ,z) \big)
\rowpl
\widetilde Z_\lambda(x).
\label{46}
\eea
where $\Theta (\lambda) \in [0, 2 \pi)$ and $\widetilde Z_\lambda \in V_{-1}^1(\Omega)$. 

\BER
\label{R3}
The singular functions \ef{20} and \ef{21} can be interpreted as "waves" travelling along 
the axis of the cusp, cf. \cite{na493} for a physical argument in a similar geometric
situation. Although such an interpretation is not directly needed in our paper, it is 
convenient to use the corresponding physical terminology, namely to call 
solutions in ker\,$ \cT_{-1}(\lambda)$ "trapped modes" and to consider $e^{i \Theta(\lambda)}$
as the "scattering coefficient" in the "diffraction" solution \ef{46}. 
\ \ $\boxtimes$
\ENR

All functions $u \in {\rm ker}\,  \cT_1(\lambda) \subset V_1^1(\Omega) \subset L^2(\Omega)$ 
belong to the domain \ef{32} of $ \cA^*$, because the inclusion $\Delta u  = - \lambda u \in 
L^2(\Omega)$  really occurs. Hence, a trapped mode is an eigenvector  corresponding to its 
eigenvalue $\lambda$ for  every self-adjoint operator $ \cA_\theta$ 
of Theorem \ref{T2}. Furthermore, it can  readily be seen that in the case
$\theta = \Theta(\lambda)$ there appears a second eigenvector \ef{46} of $ \cA_\theta$. 

\subsection{Examples of trapped modes.}
\label{sec3.2}
Following the ideas of \cite{EvLeVa} and \cite{na401} we assume that the domain 
$\Omega$ is mirror symmetric with respect to the plane $\{ x_1 = 0\}$, i.e.,
\bea
\Omega = \{ x \, : \, (-x_1, x_2, \ldots ,x_n) \in \Omega \} \label{47N}
\eea
and restrict the problem \ef{11}--\ef{12} with $f=0$ to the half 
$\Omega_+ = \{ x \in \Omega \, : \, x_1 > 0\}$ of the domain \eqref{47N},
\bea
 - \Delta u_+ (x) &=& \lambda_+ u_+(x) , \ x \in \Omega_+\ ,
\\
 \partial_\nu  u_+ (x) &=& a u_+(x) , \ x \in (\partial \Omega)_+=\partial\Omega_{+}\setminus\overline{\Sigma}, \label{47}
\eea
where we impose the artificial Dirichlet condition 
\bea
u_+(x) = 0 \ \ , \ \ \ x \in \Sigma, \label{48}
\eea
on the middle plane $\Sigma = \{ x \in \Omega \, : \, x_1 = 0\}$.

\BEL
\label{L3}
Assume that the function $u_+ \in H_{\rm loc}^1(\overline{\Omega_+} \smallsetminus
\cO )$ satisfies the inclusion $\nabla u_+ \in L^2(\Omega_+)$ and the Dirichlet condition \ef{48}. Then, the following
weighted inequality is valid:
\bea
\Vert r^{-2 }u_+ ; L^2(\Omega_+)\Vert + 
\Vert r^{-1 }u_+  ; L^2((\partial \Omega)_+)\Vert
\leq c \Vert \nabla u_+ ; L^2(\Omega_+)\Vert . \label{49}
\eea
\ENL

Proof. It suffices to verify \ef{49} in the cusp $\Pi_+^d = \{ x \in \Pi^d \, : \,
x_1 > 0\}$ and on the surface $\Gamma_+^d = \{ x \in \Gamma^d \, : \, x_1 >0 \}$. 
To this end, we write lower-dimensional inequalities in the half-section
$\omega_+^\zeta = \{ y \in \omega^\zeta \, : \, y_1 > 0 \}$,
\bea
& & \zeta^{-4} \int\limits_{\omega_+^\zeta}  |u_+(y,\zeta)|^2 dy 
\leq c \int\limits_{\omega_+^\zeta}  \big| \nabla_y  u_+(y,\zeta)\big|^2 dy  \nonumber \\
& & \zeta^{-2} \int\limits_{\partial \omega_+^\zeta}  |u_+(y,\zeta)|^2 ds_y 
\leq c \int\limits_{\omega_+^\zeta}  \big| \nabla_y  u_+(y,\zeta)\big|^2 dy 
\label{50}
\eea
coming from the Dirichlet condition \ef {48} and the 
coordinate dilatation $y \mapsto \eta = \zeta^{-2}y$. The proof is completed by
integrating \ef{50} in $\zeta \in (0,d)$ and taking into account that $ds_x = 
\big( 1 + 4z^2 |\eta \cdot \nu' (\eta)|^2\big)^{1/2} ds_y dz$. \ \ $\boxtimes$

\bigskip

The variational formulation of the problem \ef{47}--\ef{48},
\bea
(\nabla_x u_+ , \nabla_x v_+)_{\Omega_+} - a( u_+ ,  v_+)_{(\partial \Omega)_+}
= \lambda ( u_+ ,  v_+)_{\Omega_+} \ \ \forall \, v_+ \in 
H_0^1(\Omega_+;\Sigma) \label{51}
\eea
is posed in the Sobolev space $H_0^1(\Omega_+;\Sigma) \subset H^1(\Omega)$ of functions vanishing on $\Sigma$. 
Since the weight $r^{-1}$ in the second norm of \ef{49} is large when $x \to \cO$, the embedding $H^{1}_{0}(\Omega_{+},\Sigma)\subset L^{2}(\Omega_{+})$ is compact,
 and therefore the whole left-hand side of 
\ef{51} is lower semi-bounded. We deduce that the spectrum of the variational problem
\ef{51} or the differential problem \ef{47}--\ef{48} is discrete and consists of the 
monotone unbounded sequence of normal eigenvalues
\bea
\lambda_+^1 < \lambda_+^2 \leq \lambda_+^3 \leq \ldots \leq 
\lambda_+^m \leq \ldots \to + \infty. \label{52}
\eea
The corresponding eigenfunctions $u_+^m \in H_0^1(\Omega_+,\Sigma)$, $m=1,2,, \ldots$,
have smooth, odd extensions over the Dirichlet surface $\Sigma$ to $\Omega$. The extensions belong
to $H^1(\Omega)$ and thus become eigenfunctions of the original problem \ef{6},
due to the symmetry of the domain \ef{47N}. Using \ef{49} we conclude that 
the extended eigenfunctions $u^m$ belong to $V_{-2}^0 (\omega)$ and therefore
to $V_{-1}^1(\Omega)$, by Theorem 2.6 of \cite{na549}. In this way the eigenvalues 
\ef{52} are embedded into the residual spectrum of the operator $ \cA$. They also belong to
the point spectrum of any self-adjoint extension $ \cA_\theta$ of Theorem \ref{T2}.

It would be possible to show, using a general result of \cite{MaPl2,KozMaz},  
that the eigenfunctions of the problem \ef{6} belonging to $V_{-1}^1(\omega)$ 
decay exponentially as $O(e^{-\delta / r})$ for some  $\delta > 0$, although  
this argument would require  lengthy   additional computations  contained in 
\cite[Ch.\,10]{KoMaRo1} and  \cite[Sect.1.5]{na549}. However, Theorem 2.6 of 
\cite{na549}, the proof of which is much simpler,  directly implies that a solution 
$u $ of \ef{6}  belongs to $\in  V_{-\beta}^1(\Omega)$ 
for all $\beta \in \bbR_+$, i.e., it gets a super-power decay rate,  which 
will be  sufficient for our purposes.

\subsection{A peculiar property of the scattering coefficient.}
\label{sec3.3}
Let us consider the solution $Z_{\lambda + \delta}$   
of \ef{6}, when $\lambda$ is replaced by the perturbation $\lambda +\delta$, 
where $\delta $ is a small parameter with an 
arbitrary sign. We accept the simplest asymptotic ans\"atze with respect to small $\delta$
for this solution  and its scattering coefficient:
\bea
Z_{\lambda + \delta} (x) &=& Z_{\lambda} (x) + \delta Z_{\lambda}' (x)
+ \widehat  Z_{\lambda + \delta} (x) , \label{53} \\
\Theta (\lambda + \delta) &=& \Theta (\lambda) + \delta \Theta'(\lambda) + \widehat \Theta (\lambda + \delta) . \nonumber 
\eea
Both  $Z_{\lambda +\delta}$  and $\Theta (\lambda 
+\delta)$ depend smoothly on $\delta$, so that we only need to compute the correction terms, while the 
estimates  of the remainders are evident due to the general perturbation theory,
cf. \cite{Kato}. 

We derive the following problem for the function $Z_\lambda'$ 
by inserting \ef{53} to \ef{6} with $\lambda \mapsto \lambda + \delta$ and extracting terms of order
$O(\delta)$:
\bea
- \Delta Z_\lambda'(x) - \lambda Z_\lambda' (x) &=& Z_\lambda(x) \ \ , \ \ x \in \Omega, 
\nonumber \\ 
\partial_\nu  Z_\lambda' (x) &=& a Z_\lambda'(x) \ \ , \ \ x \in \partial \Omega
\smallsetminus \cO . \label{55}
\eea
Using the formulas \ef{46} with $\lambda$ and $\lambda \mapsto \lambda + \delta$,
and  the Taylor formula
\beas
e^{i \Theta (\lambda +\delta)} = e^{i \Theta (\lambda )} \big( 1 + i \delta 
\Theta'(\lambda) + O(\delta^2) \big) , 
\eeas
we derive the representation 
\bea
Z_\lambda' (x) = i \Theta'(\lambda) e^{i \Theta(\lambda)} \chi(x) \big( w_+(z) + W_\dagger(z^{-2}y,z)  \big) + \widetilde  Z_\lambda' (x) \label{56}
\eea
where  $\widetilde  Z_\lambda' \in V_{-1}^1(\Omega)$ and the incoming wave $w_-$ does not appear. 
The problem \ef{55} has a solution of the form \ef{56}.
Indeed, since the solution \ef{46} is originally defined up to a trapped mode in 
ker\,$ \cT_{-1}(\lambda)$, the orthogonality conditions 
\bea
(Z_\lambda, v)_\Omega = 0 \ \ \ \forall \, v \in {\rm ker}\,\cT_{-1}(\lambda) \label{57}
\eea
can be satisfied, and because of the Fredholm property of $ \cT_{\pm 1 } (\lambda)$, they 
guarantee the existence of a solution of the problem \ef{55} belonging to $V_1^1 (\Omega)$;
the solution is defined up to an addentum in ker\,$ \cT_1(\lambda)$, in particular, up to
the term $c_\lambda Z_\lambda$. By Theorem \ref{T1}, this general solution has the
representation \ef{29} with $w = b_+(\lambda) w_+ + b_-(\lambda)w_-$,
see \ef{19}. Setting $c_\lambda = - b_-(\lambda)$ gives \ef{56} and the orthogonality 
conditions \ef{57}. 

We now insert the functions $Z_\lambda'$ and $Z_\lambda$ into Green's formula 
in $\Omega^\zeta$ and take the limit $\zeta \to +0$ as in \ef{32} and \ef{34}.
The inhomogeneous equation in \ef{55} implies
\beas
& & \Vert Z_\lambda ; L^2(\Omega) \Vert^2 = \lim\limits_{\zeta \to +0}
\Vert Z_\lambda ; L^2(\Omega^\zeta) \Vert^2
\roweq
- \lim\limits_{\zeta \to +0} 
\int\limits_{\omega^\zeta} \big( \overline{Z_\lambda(y,\zeta)}
\partial_z Z_\lambda'(y,\zeta) - Z_\lambda'(y,\zeta) 
\partial_z \overline{ Z_\lambda(y,\zeta)} \big) dy 
\roweq
i\big( -1 \cdot 0 + \overline{e^{i \Theta(\lambda)}} i \Theta' e^{i \Theta(\lambda)}
\big) = -\Theta'. 
\eeas
Let us formulate the result.

\BET
\label{T3}
There holds the formula 
\bea
\partial_\lambda \Theta(\lambda) = -\Vert Z_\lambda ; L^2(\Omega) \Vert^2 < 0 ,
\label{58}
\eea
where $Z_\lambda$ is the solution \ef{46} of the problem \ef{6}, with \eqref{7}, subject to the
subject to the orthogonality conditions \ef{57}. 
\ENT

The relation \ef{58} means that the growth of the spectral parameter $\lambda$ in 
\ef{56} makes the scattering coefficient $e^{i \Theta(\lambda)}$ run along the 
unit circle $\{z \in  \bbC \, : \, |z|=1\} $ counter-clockwise. Unfortunately, our calculation 
does not allow to control the running speed, because the solution $Z_\lambda$
is normalized by the unit coefficient of the incoming wave $w_+$, while the norm
$\Vert Z_\lambda;L^2(\Omega)\Vert$ depends on the function
$Z_\lambda$ in the whole domain.

\section{Asymptotics of eigenvalues in the domain with a blunted cusp.}
\label{sec4}
\subsection{Formal procedure.}
\label{sec4.1}
We consider a solution of the limit problem \ef{6}, with  \eqref{7},  of the form
\ef{29} and use the coefficients $b_\pm$ in \ef{19} to satisfy the boundary condition
\ef{2} at the end $\omega^\varepsilon$ of the blunted cusp \ef{3}. We denote here by dots
terms which are inessential for our present asymptotic analysis and postpone their estimates 
to the next sections. 

In the case $a > a_\dagger$ we apply formulas \ef{20}, \ef{16} and obtain
\bea
\partial_{\nu} u(y,\varepsilon) + au(y,\varepsilon) &= &
w_0\Big(+i\mu_0 - n + \frac32 + \varepsilon a\Big) b_+ \varepsilon^{ +i\mu_0 - n +1/2}
\rowpl
w_0\Big(-i\mu_0 - n + \frac32 + \varepsilon a\Big) b_- \varepsilon^{ -i\mu_0 - n + 1/2}
+ \ldots .\label{61}
\eea
Hence, the main asymptotic term in \ef{61} vanishes provided
\bea
b_+ = B(\varepsilon)b_- \ \ , \ \ B(\varepsilon) =
\frac{-i\mu_0+ \varepsilon a - n + 3/2 }{+i\mu_0+ \varepsilon a - n + 3/2}
\varepsilon^{-2 i \mu_0}. \label{62}
\eea
The coefficient $B(\varepsilon)$ is unimodular and thus
\bea
B(\varepsilon) = e^{i T(\varepsilon)} \ \ , \ \ 
T(\varepsilon) = T_0(\varepsilon)- 2 \mu_0 \ln \varepsilon , \label{63}
\eea
where $T_0(\varepsilon)$ is a smooth function of $\varepsilon$ and $\mu_0$
is as in \ef{20} (notice that $\mu_0=0$ for $a = a_\dagger$). Comparing \ef{62}, \ef{63}
with \ef{35}, we see that the self-adjoint extension $ \cA_{\theta_\bullet^\varepsilon}$
of Theorem \ref{T2}, with the parameter
\bea
\theta_\bullet^\varepsilon = T(\varepsilon) = T_0(\varepsilon) - 2 \mu_0 \ln \varepsilon
\ \ \ (\mbox{mod}\ 2 \pi) \label{64}
\eea
is the first candidate to model the problem \ef{1}--\ef{2}.

In the case $a = a_\dagger$ we use formula \ef{21} and arrive at the relation
\bea
\partial_z u(y,\varepsilon) + a u(y,\varepsilon)&=& 
w_0\Big( \Big( - n +\frac32 + \varepsilon a \Big) (\ln \varepsilon -i)
+ 1 \Big) \varepsilon^{-n +1/2} b_+ 
\rowpl
w_0\Big( \Big( - n +\frac32 + \varepsilon a \Big) (\ln \varepsilon +i)
+ 1 \Big) \varepsilon^{-n +1/2} b_- + \ldots .  \label{65}
\eea
Deleting the main asymptotic term in \ef{65}  gives the relation
\beas
b_- = B_\dagger (\varepsilon)b_+ \ \ , \ \ B_\dagger(\varepsilon) =
\frac{+i+ \ln \varepsilon - 2(2n -3- \varepsilon  a)^{-1} }{
-i+ \ln \varepsilon - 2(2n -3- \varepsilon  a)^{-1}}
=: e^{i T_\dagger(\varepsilon)}. 
\eeas
The corresponding parameter 
\bea
\theta_\bullet^\varepsilon = T_\dagger (\varepsilon)  \label{67}
\eea
of the appropriate self-adjoint extension $ \cA_{\theta_\bullet^\varepsilon}$
for modelling the problem \ef{1}--\ef{2} behaves in a very different way as $\varepsilon \to
+0$ in comparison with the function \ef{64}, which is "almost linear" in $\ln \varepsilon$,
namely we have
\beas
B_\dagger (\varepsilon) = 1 + O( |\ln \varepsilon|^{-1} ) \ \ , \ \ 
T_\dagger = O( |\ln \varepsilon|^{-1} ) \ \ \ \mbox{as} \ \varepsilon \to +0.
\eeas

\subsection{Operator formulation of the problem.}
\label{sec4.2}
We consider the integral identity  \cite{Lad}
\bea
( \nabla u^\varepsilon , \nabla v^\varepsilon)_{\Omega^\varepsilon}
- a( u^\varepsilon ,  v^\varepsilon)_{\Omega^\varepsilon} = \lambda^\varepsilon
( u^\varepsilon , v^\varepsilon)_{\Omega^\varepsilon} \ \ 
\forall \, v^\varepsilon \in H^1(\Omega^\varepsilon) \label{69}
\eea
for the problem \ef{1}--\ef{2}; notice that $\Omega^\varepsilon$ is a Lipschitz 
domain.

The following inequality, where $c$ is independent of $u^\varepsilon \in 
H^1(\Omega^\varepsilon)$, can be verified along the same lines as Proposition in \cite{na401}:
\beas
\Vert r^{-1} u^\varepsilon;L^2(\Omega^\varepsilon) \Vert \leq c
\Vert u^\varepsilon; H^1(\Omega^\varepsilon) \Vert. 
\eeas
Thus, the standard norm of $H^1(\Omega^\varepsilon)$ is equivalent with the weighted norm
$\Vert \cdot ; V_0^1(\Omega^\varepsilon)\Vert$, uniformly in $\varepsilon$, see \ef{23}.

We need some estimates in order to write an abstract formulation of \ef{69}.

\BEL
\label{L4}
The trace inequalities
\bea
\Vert u^\varepsilon; L^2(\partial\Omega^\varepsilon \smallsetminus \omega^\varepsilon) \Vert^2
&\leq& c\big( \delta\Vert  \nabla u^\varepsilon; L^2(\Omega^\varepsilon) \Vert^2 
+ (1 + \delta^{-1} ) \Vert r^{-1} u^\varepsilon; L^2(\Omega^\varepsilon) \Vert^2 \big),
\label{71} \\
\Vert u^\varepsilon; L^2(\omega^\varepsilon) \Vert &\leq &
c \sqrt{\varepsilon}\Vert u^\varepsilon; V_0^1(\Omega^\varepsilon) \Vert  \label{72}
\eea
hold true with  constants $c$ depending  on neither $u^\varepsilon \in 
H^1(\Omega^\varepsilon)$ nor $\varepsilon \in (0, \varepsilon_0)$. 
\ENL

Proof. The inequality \eqref{72} is verified in \cite{NaTaNEW}, Lemma 5.1.
Concerning \eqref{71}, it is enough to check the statement for smooth real-valued functions and by replacing
$\Omega^\varepsilon \mapsto \Pi^d \smallsetminus \Pi^\varepsilon$ and
$\partial \Omega^\varepsilon \mapsto \Gamma^d \smallsetminus \Gamma^\varepsilon$.  
To this end we make the coordinate compression $\eta \mapsto y = z^2 \eta$, cf. \ef{4},
in the standard trace inequality
\beas
\Vert U ; L^2(\partial \omega) \Vert^2 \leq c_\omega 
\Vert U ; H^1 ( \omega) \Vert \, \Vert U ; L^2(\omega) \Vert
\eeas
and obtain
\bea
& & \Vert u^\varepsilon ; L^2(\partial \omega^z) \Vert^2 
\rowleq 
c_\omega z^{-2} \big( z^2\Vert \nabla_y u^\varepsilon ; L^2(\omega^z) \Vert^2
+ \Vert u^\varepsilon ; L^2(\omega^z) \Vert^2 \big)^{1/2}
\Vert u^\varepsilon ; L^2(\omega^z) \Vert 
\rowleq
c_\omega \big( z \Vert \nabla_y u^\varepsilon ; L^2(\omega^z) \Vert
+ z^{-1} \Vert u^\varepsilon ; L^2(\omega^z) \Vert \big)
z^{-1} \Vert u^\varepsilon ; L^2(\omega^z) \Vert
\rowleq
C_\omega \Big( \delta \Vert \nabla_y  u^\varepsilon ; L^2(\omega^z) \Vert^2
+ \Big( 1 + \frac1\delta \Big) \Vert r^{-1} u^\varepsilon ; L^2(\omega^z) \Vert^2
\Big)  \label{73}
\eea
Here we replaced $z$ with $d$ in front of $\nabla_y u^\varepsilon$, inserted 
$z^{-1}$ as $r^{-1}$ inside the Lebesgue norms of $u^\varepsilon$, and applied the 
Cauchy inequality $2ab \leq \delta a^2 + \delta^{-1}b^2$. Finally, \ef{71} follows by 
integrating  in $z \in (\varepsilon, d)$ the inequality obtained in \ef{73}.
 \ \ $\boxtimes$

\bigskip

We introduce in the Hilbert space $\cH^\varepsilon =H^1(\Omega^\varepsilon)$ the new
scalar product
\bea
\langle u^\varepsilon, v^\varepsilon \rangle_\varepsilon = 
( \nabla u^\varepsilon, \nabla v^\varepsilon)_{\Omega_\varepsilon}
-a( u^\varepsilon,  v^\varepsilon)_{\partial\Omega_\varepsilon}
+ \ell \varepsilon^{-2} ( u^\varepsilon,  v^\varepsilon)_{\Omega_\varepsilon};
\label{75}
\eea
the properties of a scalar product follow for a large enough $\ell >0$
from Lemma \ref{L4} and the obvious relation $ r > \varepsilon$ in 
$\Omega^\varepsilon$. Moreover, we can and do fix $\ell $ such that 
\bea
\Vert u^\varepsilon; \cH^\varepsilon \Vert^2 =
\langle u^\varepsilon, u^\varepsilon \rangle_\varepsilon
\geq \varepsilon^{-2} \Vert u^\varepsilon; L^2(\Omega^\varepsilon) \Vert^2.
\label{75N}
\eea
We define the operator $\cK^\varepsilon$ in $\cH^\varepsilon$ by the identity
\bea
\langle \cK^\varepsilon u^\varepsilon, v^\varepsilon \rangle_\varepsilon
=(u^\varepsilon, v^\varepsilon )_{\Omega^\varepsilon}
\ \ \forall\, u^\varepsilon, v^\varepsilon \in \cH^\varepsilon ,
\label{76}
\eea
so that $\cK^\varepsilon$ becomes positive, continuous and symmetric, hence, self-adjoint.
Moreover, it is compact, and by \cite[Thm.\,10.1.5.]{BiSo}, its essential spectrum 
coincides with $\{\mu=0\}$, and according to \ef{75} and \ef{76}, problem \ef{69} is 
equivalent with the abstract equation
\beas
\cK^\varepsilon u^\varepsilon = \sfk^\varepsilon u^\varepsilon \ \ \ 
\mbox{in} \ \cH^\varepsilon   
\eeas
with the new spectral parameter
\bea
\sfk^\varepsilon = (\ell \varepsilon^{-2} + \lambda^\varepsilon)^{-1}
\label{78}
\eea
The discrete spectrum of $\cK^\varepsilon$,
\bea
\sfk_1^\varepsilon \geq \sfk_2^\varepsilon \geq 
\ldots \geq \sfk_m^\varepsilon \geq \ldots \to +0 , \label{79}
\eea
is related to the eigenvalue sequence \ef{5} via formula \ef{78}.

The next assertion is known as the lemma on "near eigenvalues", \cite{ViLu},
and it follows from the spectral decomposition of a resolvent,
see \cite[Ch.\,6]{BiSo}.

\BEL
\label{LL4}
Let $u_\bullet^\varepsilon \in H^1(\Omega^\varepsilon)$, 
$u_\bullet^\varepsilon \not=0$, 
and $\sfk_\bullet^\varepsilon \in \bbR$ satisfy
\bea
\Vert \cK^\varepsilon u_\bullet^\varepsilon - \sfk_\bullet^\varepsilon
u_\bullet^\varepsilon ; \cH^\varepsilon\Vert = 
\delta_\bullet^\varepsilon \Vert u_\bullet^\varepsilon ; \cH^\varepsilon\Vert 
\ , \ \ \delta_\bullet^\varepsilon \in [0, k_\bullet^\varepsilon).
\label{80}
\eea
Then,  there exists an eigenvalue $\sfk^\varepsilon_m$ belonging to the sequence \ef{79}
such that
\beas
| \sfk_\bullet^\varepsilon - \sfk_m^\varepsilon | \leq \delta_\bullet^\varepsilon.
\eeas
\ENL

\subsection{Error estimates for the approximation of the spectrum in $\Omega^\varepsilon$.}
\label{sec4.3}
We fix a spectral parameter $\lambda\in \bbR$, assume \eqref{7} and 
consider the solution \ef{46} of the problem \ef{6} with the scattering 
coefficient  $e^{i \Theta(\lambda)}$, and find a positive sequence $\{ \varepsilon_k 
\}_{k \in \bbN}$ tending to 0 such that 
\bea
\Theta(\lambda) = T_0(\varepsilon_k) - 2\mu_0 \ln \varepsilon_k \ \ \ 
(\mbox{mod} 2 \pi). \label{82}
\eea
Then, $u_\bullet^{\varepsilon_k} = Z_\lambda$ is an eigenfunction 
corresponding to the eigenvalue $\lambda$  of the self-adjoint extension $ 
\cA_{\theta_\bullet^{\varepsilon_k}}$ of Theorem \ref{T2}. Notice that the 
exponent $\Theta(\lambda)$ depends continuously on the parameter $\lambda$ (to 
see this recall that the space $\bfV_\pm^1(\Omega)
\supset \cD( \cA^*)$ can be identified with the direct product \ef{c2v} and 
apply general results on non-selfdajoint linear operators in \cite[Ch.\,1]
{GoKr},  \cite[Ch.\,4]{Kato}), and therefore
\beas
\Vert Z_\lambda ; \bfV_{\pm 1}^1(\Omega)\Vert \leq C_Z \ \ , \ \
\Vert Z_\lambda ; L^2(\Omega)\Vert \geq c_Z > 0  
\eeas
uniformly with respect to $\lambda$ in a compact set.  Hence, in view of \ef{75N}, we have
\bea
\Vert u_\bullet^{\varepsilon_k} ; \cH^\varepsilon \Vert^2 \geq 
\varepsilon_k^{-2} \Vert u_\bullet^{\varepsilon_k} ; 
L^2(\Omega^\varepsilon)\Vert^2 
\geq \frac12 c_Z^2 \varepsilon_k^{-2}  \label{83N}
\eea
because
\beas
& & \Vert u_\bullet^{\varepsilon_k} ; L^2(\Pi^\varepsilon ) \Vert^2 
\leq c \varepsilon_k^{2(1- \gamma)} 
\Vert z^{\gamma-1} u_\bullet^{\varepsilon_k} ; L^2(\Pi^\varepsilon ) \Vert^2  
\rowleq 
c\Vert \varepsilon_k^{2(1- \gamma)}  u_\bullet^{\varepsilon_k} ; 
V_\gamma^1 (\Omega  ) \Vert^2
\leq c \varepsilon_k^{2(1- \gamma)} 
\Vert Z_\lambda ; \bfV_{\pm 1}^1(\Omega)\Vert
\eeas
where $\gamma \in (0,1)$ and $\cD(\cA_{\theta_\bullet^{\varepsilon_k}} )
\subset \cD(\cA^*) \subset V_\gamma^1(\Omega)$, cf. Remark \ref{R1}. 

We are going to prove that the problem \eqref{1}--\eqref{2} with $a > a_\dagger$ in the domain $\Omega^{\varepsilon_k}$
has an eigenvalue $\lambda^{\varepsilon_k}$ in the vicinity of $\lambda$.
In what follows we write $\varepsilon$ instead of $\varepsilon_k$. The threshold case $a = a_\dagger$ as well as other eigenvalues with stable asymptotics will be considered in Section 
\ref{sec4.4}. 

According to \ef{78} we set
\bea
\sfk_\bullet^\varepsilon = ( \ell \varepsilon^{-2} + \lambda)^{-1} . \label{83}
\eea
To compute the factor $\delta_\bullet^\varepsilon$ in \ef{80}, we use the 
definition of the norm of a  Hilbert space, \ef{75} and \ef{76}, and write
\bea
& & \Vert \cK^\varepsilon u_\bullet^\varepsilon - \sfk_\bullet^\varepsilon 
u_\bullet^\varepsilon ; \cH^\varepsilon \Vert = 
\sup \big| \langle \cK^\varepsilon u_\bullet^\varepsilon - \sfk_\bullet^\varepsilon 
u_\bullet^\varepsilon ,v^\varepsilon \rangle_\varepsilon \big|
\roweq
\varepsilon^2 (\ell + \varepsilon^2 \lambda)^{-1} 
\sup | (\ell \varepsilon^{-2} + \lambda
) (u_\bullet^\varepsilon, v^\varepsilon)_{ \Omega^\varepsilon}
- (\nabla u_\bullet^\varepsilon, \nabla v^\varepsilon)_{ \Omega^\varepsilon}
\rowpl
a(u_\bullet^\varepsilon, v^\varepsilon)_{\partial \Omega^\varepsilon}
- \ell \varepsilon^{-2} (u_\bullet^\varepsilon, v^\varepsilon)_{ \Omega^\varepsilon} \big|
\roweq
\varepsilon^2 (\ell + \varepsilon^2 \lambda)^{-1}
\sup \big| ( \Delta u_\bullet^\varepsilon + \lambda  u_\bullet^\varepsilon, v^\varepsilon)_{ \Omega^\varepsilon}
+ (\partial_n u_\bullet^\varepsilon - a  u_\bullet^\varepsilon, 
v^\varepsilon)_{\partial  \Omega^\varepsilon} \big| . \label{84}
\eea
Here, the supremum is taken over the unit ball $\{ v^\varepsilon 
\in \cH^\varepsilon \, : \, \Vert v^\varepsilon ; \cH^\varepsilon \Vert \leq 1  \}$
and the Green formula was applied. Furthermore, on the last line the scalar 
product in $L^2(\Omega^\varepsilon)$ is null, and also 
\beas
(\partial_n u_\bullet^\varepsilon - a  u_\bullet^\varepsilon, 
v^\varepsilon)_{\partial  \Omega^\varepsilon 
\smallsetminus \omega^\varepsilon} = 0.
\eeas
So we are left with only
\bea
I_\bullet^\varepsilon(v^\varepsilon) = - \int\limits_{\omega^\varepsilon}
\overline{v^\varepsilon(y,\varepsilon)} \big( \partial_z u_\bullet^\varepsilon 
(y,\varepsilon) + a u_\bullet^\varepsilon(y,\varepsilon) \big) dy. \label{85}
\eea
According to our preparatory calculation \ef{61} of the parameter $\theta_\bullet^\varepsilon$ 
in \ef{64}, the main asymptotic term $w_- + e^{i \Theta (\lambda)}
w_+$ disappears in the integrand in \ef{85} so that the integral itself
reduces to the sum
\beas
I_W^\varepsilon(v^\varepsilon) + \widetilde I^\varepsilon(v^\varepsilon) 
& =&  
- \int\limits_{\omega^\varepsilon} 
\overline{v^\varepsilon(y,\varepsilon)} \Big( \frac{d}{dz} + a \Big) 
W(z^{-2} y ,z) \Big|_{z=\varepsilon} dy 
\rowmi
 \int\limits_{\omega^\varepsilon} 
\overline{v^\varepsilon(y,\varepsilon)} 
\big( \partial_z \widetilde u_\bullet^\varepsilon 
(y,\varepsilon) + a \widetilde u_\bullet^\varepsilon(y,\varepsilon) \big) dy.
\eeas
Recalling the form \ef{14}  of the correction terms $W_-(z^{-2} y ,z ) +
e^{i \Theta( \lambda)} W_+(z^{-2} y ,z )$ we get
\beas
\frac{d}{dz} W(z^{-2} y ,z) &=&-\frac{2}{z^3}y\cdot\nabla_{\eta}W(z^{-2} 
y ,z)+\frac{\partial W}{\partial z} \left(z^{-2}y,z\right)
\\
&=&-2z^{\mu-1}y\cdot\nabla_{\eta}W_{0}(z^{-2} y)+z^{\mu+1}W_{0}(z^{-2}y)
\eeas
For $y\in \omega^{\varepsilon}$ there holds $|y| \leq c \varepsilon^{2}$, therefore, 
$$\Big|\Big( \frac{d}{dz} + a \Big)  W(z^{-2} y ,z) \Big|_{z=\varepsilon}\Big| \leq c \varepsilon^{1-n+3/2}(|W_{0}\big(\frac{y}{\varepsilon^2}\big)|+|\nabla_{\eta}W_{0}\big(\frac{y}{\varepsilon^2}\big)|) $$
Integrating over $y\in \omega^{\varepsilon}$ yields
$$\Big\| \Big( \frac{d}{dz} + a \Big)  W(z^{-2} y ,z) \Big|_{z=\varepsilon};L^{2}(\omega_{\varepsilon}) \Big\| \leq c\varepsilon^{1-n+\frac{3}{2}}\varepsilon^{n-1}\|W_{0};H^{1}(\omega)\|=\varepsilon^{3/2}\|W_{0};H^{1}(\omega)\|$$
The inequality \ef{72} implies 
\beas
 |I_W^\varepsilon(v^\varepsilon)| 
\leq 
c\varepsilon^{1/2} \varepsilon^{3/2} =c \varepsilon^2.
\eeas
Since $\widetilde u_\bullet \in V_{-1}^1(\Omega)$, we again apply \ef{72} to obtain
\bea
& & \Big| \int\limits_{\omega^\varepsilon} \overline{v^\varepsilon(y,\varepsilon)}
\widetilde u_\bullet^\varepsilon(y,\varepsilon) dy \Big| 
\leq  \Vert v^\varepsilon; L^2(\omega^\varepsilon)\Vert \,
\Vert \widetilde u_\bullet^\varepsilon; L^2(\omega^\varepsilon)\Vert
\rowleq
c \varepsilon^{1/2} \Vert v^\varepsilon; \cH^\varepsilon \Vert \varepsilon
\Vert r^{-1}  \widetilde u_\bullet^\varepsilon; L^2(\omega^\varepsilon)\Vert 
\rowleq 
c\varepsilon^2 \Vert r^{-1}  \widetilde u_\bullet^\varepsilon; V_0^1(\Omega^\varepsilon)\Vert 
\leq c\varepsilon^2 \Vert \widetilde u_\bullet^\varepsilon; V_{-1}^1(\Omega^\varepsilon)\Vert . \label{87}
\eea
However, the estimation of the integral with $\partial_{z}\widetilde u_\bullet^\varepsilon$
is much more involved. Indeed, a direct application of the weighted trace
inequality in Lemma \ref{L4}  does not help, because \cite[Lem.\,3.2]{na549} 
proves that 
the second derivatives $\nabla^2 \widetilde u_\bullet^\varepsilon$ belong 
to $V_1^0(\Omega)$, not to $V_0^0(\Omega)$ as in the Kondratiev theory, cf. also
Remark \ref{R2}. However, the situation can be improved in the case of the trivial
right-hand side $f=0$. Indeed, in our case the function $u_\bullet^\varepsilon$
satisfies the homogeneous problem \ef{1}--\ef{2} so that the one can use the procedure
in \cite[Sec.\,2]{na549} to construct several higher-order asymptotic terms
\bea
W^{(k)} ( z^{-2}y,z) = z^{2k -n + i \mu_0 +3/2} W_ {0+}^{(k)} (z^{-2}y)  
+  z^{2k -n  -i \mu_0 +3/2} W_ {0-}^{(k)}  (z^{-2}y)  \label{88}
\eea
(notice that $W^{(0)}$ coincides with $W$) and to obtain a remainder $\widetilde u_\bullet^{\varepsilon(k)}$ with any 
desired power-law decay rate as $z \to +0$. In particular, $\partial_z 
\widetilde u_\bullet^{\varepsilon(1)} \in V_{-1}^1(\Omega)$. The additional 
term
$W^{(1)}$ can be processed similarly to the calculation \ef{6}, while the integral
$$
\int\limits_{\omega^\varepsilon} \overline{v^\varepsilon(y,\varepsilon)}
\partial_z \widetilde u_\bullet^{\varepsilon(1)} (y, \varepsilon)dy
$$ 
can be estimated in the same way as in \ef{87}. Summing up, we conclude that the
order of the norm \ef{84} with respect to  $\varepsilon$ is determined by the 
main correction term \ef{88} at $k=0$. Thus, collecting the estimates derived so far 
and comparing \ef{80} with \ef{83N} we conclude that
\bea
\delta_\bullet^\varepsilon \leq c \varepsilon^5 (\ell + \varepsilon^2 \lambda)^{-1}.
\label{89}
\eea

\BET
\label{AS1}
Assume $a > a_\dagger$. Let $\lambda\in \bbR$ be fixed, and let the sequence $\{ \varepsilon_k\}_{k \in\bbN}$ satisfiy \ef{82}. Then, there exist $\varepsilon_\bullet >0$ and   $c_\bullet>0$ such that for $\varepsilon_k \leq 
\varepsilon_\bullet$ the problem \ef{1}--\ef{2} with $\varepsilon=\varepsilon_k$
has an eigenvalue $\lambda_{m(\varepsilon_k)}^{\varepsilon_k}$, \ef{5},
satisfying the estimate
\beas
\big| \lambda_{m(\varepsilon_k)}^{\varepsilon_k} - \lambda \big| \leq c_\bullet
\varepsilon_k .  
\eeas
\ENT

Proof. According to Lemma \ref{LL4} and formula \ef{89} we find an
eigenvalue $\sfk_{m(\varepsilon_k)}^{\varepsilon_k} $ of the operator
$\cK^{\varepsilon_k}$, see \ef{79} and \ef{76}, such that
\beas
\big| \sfk_{m(\varepsilon_k)}^{\varepsilon_k} - \sfk_\bullet^{\varepsilon}\big| \leq 
c\varepsilon^5 (\ell + \varepsilon^2 \lambda)^{-1}.
\eeas
Taking into account \ef{78} and \ef{83} we can write
\bea
\big|  \lambda_{m(\varepsilon_k)}^{\varepsilon_k} - \lambda \big|
&\leq & c \varepsilon_k^5 (\ell +\varepsilon_k^2\lambda)^{-1} 
(\ell \varepsilon_k^{-2} + \lambda)
\big(\ell \varepsilon_k^{-2} + \lambda_{m(\varepsilon_k)}^{\varepsilon_k} \big)
\roweq
c \varepsilon_k^3 \big(\ell \varepsilon_k^{-2} + 
\lambda_{m(\varepsilon_k)}^{\varepsilon_k} \big) .
\label{91}
\eea
Fixing a small enough $\varepsilon_\bullet \leq 1$, namely $c \varepsilon_\bullet^3
\leq 1/2$, we obtain from \ef{91} that
\beas
\lambda_{m(\varepsilon_k)}^{\varepsilon_k}  \leq \lambda + c \varepsilon_k
\big(\ell + \varepsilon_k^2 \lambda_{m(\varepsilon_k)}^{\varepsilon_k}  \big) 
\ \ \Rightarrow \ \ \lambda_{m(\varepsilon_k)}^{\varepsilon_k}  \leq 
2\lambda+ \ell \leq 
2 |\lambda| + \ell. \ \  \boxtimes
\eeas

\subsection{Threshold case and "stable" eigenvalues.}
\label{sec4.4}

At $a = a_\dagger$ 
the spectral parameter $\lambda\in \bbR$ still gives rise to 
the exponent $\Theta(\lambda)$ of the scattering coefficient in \ef{46}, but now the sequence $\{ \varepsilon_k\}_{k \in\bbN}$ will be defined by
\bea
\Theta(\lambda) = T_\dagger (\varepsilon)  \label{67bis}.
\eea
The waves \ef{21} include the logarithmic factor $\ln z$, and the logarithm only 
causes only self-evident technical differences in the calculations and arguments 
in Section \ref{sec4.3}. Hence, we just reformulate Theorem \ref{AS1} as follows.

\BET
\label{AS3}
Assume $a =a_\dagger$. Let $\lambda\in \bbR$ be fixed, and let the sequence $\{ \varepsilon_k\}_{k \in\bbN}$ satisfiy \ef{67bis}. Then, there exist $\varepsilon_\bullet >0$ and   $c_\bullet>0$ such that for $\varepsilon_k \leq 
\varepsilon_\bullet$ the problem \ef{1}--\ef{2} with $\varepsilon=\varepsilon_k$
has an eigenvalue $\lambda_{m_{\dagger}(\varepsilon_k)}^{\varepsilon_k}$, \ef{5},
satisfying the estimate
\beas
\big| \lambda_{m_\dagger(\varepsilon_{k})}^{\varepsilon_{k}} - \lambda \big| \leq c_\bullet
\varepsilon_{k}(1 + |\ln \varepsilon_{k}|) . 
\eeas
\ENT

If $u^{tr} \in V_{-1}^1(\Omega)$ is a trapped mode for the problem \ef{6} with
some $\lambda\in \bbR$ (Remark \ref{R3} and Section \ref{sec3.2}), then this point $\lambda$
is an eigenvalue of every self-adjoint extension $\cA_\theta$, $\theta \in [0, 2\pi)$,
see Theorem \ref{T2}. Moreover, since $u^{tr} \in V_{-\beta}^1(\Omega)$ for all weight 
indices $\beta$, repeating the calculations and arguments in the previous 
section leads to the following assertion, which includes the threshold case too, because a trapped mode has the same
fast decay properties both in the case $a=a_{\dagger}$ and $a>a_{\dagger}$. 

\BET
\label{AS2}
Assume that $\lambda \in \bbR$ is an embedded eigenvalue of the problem \eqref{6}, $a \geq a_\dagger$, related  with the trapped mode $u^{tr} \in V_{-\beta}^1(\Omega)$. Then, for any $N \in \bbR_+$ there exist 
$\varepsilon_N >0$ and $c_N> 0$ such that, for all $\varepsilon \in (0,\varepsilon_{N})$,
the problem \ef{1}--\ef{2} has an eigenvalue $\lambda_{m(\varepsilon)}^{\varepsilon}$,
\ef{5}, satisfying the estimate
\beas
\big| \lambda_{m(\varepsilon)}^{\varepsilon} - 
\lambda \big| \leq c_N \varepsilon^N . 
\eeas
\ENT

\section{Conclusions and possible generalizations.}
\label{sec5}
\subsection{Asymptotic behavior of eigenvalues above threshold.}
\label{sec5.1}
According to Theorem \ref{AS2}, any trapped mode $u^{tr} \in V_{-1}^1
(\Omega)$ of the problem \ef{6} with $\lambda = \lambda^{tr}$ gives rise to a
family $\{ \lambda_{m(\varepsilon)}^\varepsilon \}_{\varepsilon \in (0, \varepsilon_{N})}$
of eigenvalues of the problems \ef{1}--\ef{2} in $\Omega^\varepsilon$ staying
in the $c_N \varepsilon^N$-neighborhood of the point $\lambda^ {tr}$. 
We call the eigenvalue $\lambda_{m(\varepsilon)}^{\varepsilon}$ as a stable one in spite of 
the fact that the number $m(\varepsilon)$ of these eigenvalues in the sequence \ef{5} changes 
infinitely many times, when $\varepsilon \to +0$ (see the explanation in the next paragraph).

We have detected eigenvalues  in the spectra of the problems \ef{1}--\ef{2} 
having completely different behavior as $\varepsilon\to +0$. Indeed, let 
$\lambda\in \bbR$ be fixed. Theorem \ref{AS1} shows that in a neighborhood of 
$\lambda$ periodically, with the  period
\beas
\pi^{-1} m_0 
\eeas
in the logarithmic scale $|\ln \varepsilon|$, there appears an eigenvalue of
the problem \ef{1}--\ef{2} in the domain $\Omega^\varepsilon$, which grows 
because the length of the broken piece $\Pi^\varepsilon$ diminishes. This
eigenvalue crosses a neighborhood of $\lambda$ at a high speed
$O(|\ln \varepsilon|)$ (in particular, $\lambda$ becomes an eigenvalue of the 
problem \ef{1}--\ef{2}); notice that this is the very reason for the rapid 
changes of the number $m(\varepsilon)$ of the stable eigenvalues mentioned
above. In other words, any point of the real axis above the threshold
becomes  a "blinking eigenvalue", as $\varepsilon \to +0$. 

We may change the point of view and watch over eigenvalues the eigenfunctions 
 of which are of the form \ef{46} with the exponent
\bea
\Theta (\lambda) = T_0(\varepsilon) - 2 \mu_0
\ln \varepsilon = 2 \mu_0 |\ln \varepsilon| + T_0(0) + O(\varepsilon)
\label{99}
\eea
of the scattering coefficient. The function \ef{99} is monotone growing when $|\ln 
\varepsilon| \to  + \infty$. Hence, the scattering coefficient 
$e^{i \Theta(\lambda)}$ moves counter-clockwise, when $\varepsilon \to +0$.
By Theorem \ref{T3}, such movement of the coefficient corresponds to the
monotone descend of the spectral parameter down along the real axis.

\subsection{Other boundary conditions.}
We have imposed in Section \ref{sec1} the same Robin condition \eqref{2} both 
on the blunted surface of the peak $\Pi^d \smallsetminus \overline{\Pi^\varepsilon}$
and on the massive part $\partial \Omega^\varepsilon \smallsetminus 
\overline{\omega^\varepsilon} $ of the boundary. Replacing 
\eqref{2} 
by
\bea
\partial u^\varepsilon(x) &=& au^\varepsilon(x) , \ \ x \in 
\partial \Omega^\varepsilon \smallsetminus \overline{\omega^\varepsilon} ,
\label{abo1} \\
\partial u^\varepsilon(x) &=& 0 , \ \ x \in \omega^\varepsilon \label{abo2}
\eea
does not cause any changes to the calculations and justification. 
Replacing the Neumann condition \eqref{abo2} by the Dirichlet one,
\bea
u^\varepsilon(x) &=& 0 , \ \ x \in \omega^\varepsilon , \label{abo3}
\eea
the properties of the problem \eqref{1}, \eqref{abo1}, \eqref{abo3} 
still remain very similar, although our calculation of the extension 
parameters \eqref{64} and \eqref{67} requires a minor (simplifying)
modification.

\subsection{Other shapes.}
All the results on the problems \eqref{1}, \eqref{2} and \eqref{1}, \eqref{abo1}, 
\eqref{abo2} or  \eqref{abo3} remain unchanged, if  the straight end 
$\omega^{\varepsilon} =\{ x \in \Pi^d \, : \, z = \varepsilon \}$ of the domain \eqref{3} is made into a curved one, i.e. 
$$
\Omega^{\varepsilon}=\{x=(y,z): z>\varepsilon+\varepsilon^2H(\varepsilon^{-2}y)\}
$$ 
where $H$ is a Lipschitz function in $\overline\omega$.

One may also divide the lateral boundary of $\Gamma^d$ of the cusp $\Pi^d$ into 
two non-empty and non-intersecting parts $\Gamma_k^d = \{ x \, : \,
z \in (0,d), \ z^{-2} y \in \gamma_k \}$, where 
both  sets $\gamma_k$, $k=1,2$, are open submanifolds of $\partial \omega$
and $\partial \omega= \overline{ \gamma_1} \cup  \overline{ \gamma_2} $. 
If one keeps 
the Robin condition on $\Gamma_1^d$ and imposes the Neumann condition on 
$\Gamma_2^d$, the above-discovered properties of the spectrum are still retained
by the modified problem. However, in the case of the Dirichlet condition on 
$\Gamma_2^d$ the spectrum of the problem on $\Omega^\varepsilon$
is  discrete and therefore its eigenvalues are hardly movable. 
In  general,  changes of the boundary conditions outside a neighborhood of the 
tip $\cO$ do  not affect the above described properties of the spectrum.

\end{document}